\documentclass[11pt]{article}
\usepackage{amsmath}
\usepackage{amsthm}
\usepackage{epsfig}
\usepackage{leftidx}
\usepackage{graphicx,subfigure}
\usepackage{amsfonts}
\usepackage{graphicx, color, epstopdf}
\usepackage{cite}
\usepackage{flafter}
\usepackage{amssymb,amsmath}%\eqref
\usepackage{multicol}
\usepackage{booktabs}
\usepackage{multirow}
\usepackage{enumerate}
\usepackage{enumitem}
\usepackage{algpseudocode,algorithm,algorithmicx}
\newtheorem{theorem}{Theorem}[section]
\newtheorem{lemma}{Lemma}[section]

\newtheorem{remark}{Remark}
\newtheorem{definition}{Definition}

\newcommand{\rmk}{\mathrm{k}}

\newcommand{\defeq}{:=}
\newcommand{\zd}{\,\mathrm{d}}
\newcommand{\abs}[1]{\left|#1\right|}
\newcommand{\absb}[1]{\big|#1\big|}
\newcommand{\absB}[1]{\Big|#1\Big|}

\newcommand{\bra}[1]{\left(#1\right)}
\newcommand{\brab}[1]{\big(#1\big)}
\newcommand{\braB}[1]{\Big(#1\Big)}
\newcommand{\brat}[1]{(#1)}
\newcommand{\kbra}[1]{\left[#1\right]}
\newcommand{\kbrab}[1]{\big[#1\big]}
\newcommand{\kbraB}[1]{\Big[#1\Big]}

\newcommand{\mynorm}[1]{\left\|#1\right\|}
\newcommand{\mynormb}[1]{\big\|#1\big\|}

\title{Stability of variable-step BDF2 and BDF3 methods}
\author{
Zhaoyi Li\thanks{Department of Mathematics,
Nanjing University of Aeronautics and Astronautics,
Nanjing 211106, P. R. China. Email: 1525170777@qq.com}
\quad
Hong-lin Liao\thanks{Corresponding author. ORCID 0000-0003-0777-6832. Department of Mathematics,
Nanjing University of Aeronautics and Astronautics,
Nanjing 211106, China; Key Laboratory of Mathematical Modelling
and High Performance Computing of Air Vehicles (NUAA), MIIT, Nanjing 211106, China.
 Hong-lin Liao (liaohl@nuaa.edu.cn and liaohl@csrc.ac.cn)
is supported by a grant 12071216 from
National Natural Science Foundation of China.}}
\date{\today}

\begin{document}

\maketitle

\begin{abstract}
We prove that the two-step backward differentiation formula (BDF2)
method is stable on arbitrary time grids; 
while the variable-step BDF3 scheme is stable if almost all adjacent step ratios are less than 2.553.
These results relax the severe mesh restrictions in the literature and 
provide a new understanding of variable-step BDF methods.
Our main tools include the discrete orthogonal convolution kernels
and an elliptic-type matrix norm.
\\[1ex]
\\[1ex]
\end{abstract}

\section{Introduction}
\setcounter{equation}{0}

In this paper, we revisit the stability of variable-step
backward differentiation formulas (BDF) for the following initial value problem
\begin{align*}
\frac{\zd v}{\zd t}=f(t,v)\quad\text{for $0<t\le T$}.
\end{align*}
Choose time levels $0=t_0<t_1<t_2<\cdots<t_N=T$ with the variable time-step
$\tau_k\defeq t_k-t_{k-1}$ for $1\le k\le N$, and the maximum step size
$\tau\defeq\max_{1\le k\le N}\tau_k$. Put the adjacent time-step ratios
$$r_{k}\defeq\frac{\tau_{k}}{\tau_{k-1}}\quad\text{for $2\le k\le N$}\quad\text{with}\quad r_1:=0.$$
For any sequences $\{v^n\}_{n=0}^N$, denote $\partial_{\tau}v^n:=(v^n-v^{n-1})/\tau_n$.
Let $\Pi_{n,k}v$ ($k=1,2,3$) be the interpolating polynomial of a function~$v$
over the nodes $t_{n-k}$, $\cdots$, $t_{n-1}$~and $t_{n}$.
To describe the variable coefficients, define the following three functions
\begin{align}
d_0(x,y):=&\,\frac{1+2x}{1+x}+\frac{xy}{1+y+xy},\label{funs: d0}\\
d_1(x,y):=&\,-\frac{x}{1+x}
-\frac{xy}{1+y+xy}
-\frac{xy^2}{1+y+xy}\frac{1+x}{1+y},\label{funs: d1}\\
d_2(x,y):=&\,\frac{xy^2}{1+y+xy}
\frac{1+x}{1+y},\qquad\text{for $x, y\ge0$.}\label{funs: d2}
\end{align}
Taking $v^n=v(t_n)$, one has the BDF1 formula
$D_1v^n:=\bra{\Pi_{n,1}v}'(t_n)=d_0(0,0)\partial_{\tau} v^n$ for $n\ge1$,
and the variable-step BDF2 formula
\begin{align}\label{eq: BDF2 formula}
D_2v^n:=\bra{\Pi_{n,2}v}'(t_n)=d_0(r_n,0)\partial_{\tau} v^{n}
+d_1(r_n,0)\partial_{\tau} v^{n-1}\quad\text{for $n\ge2$}.
\end{align}
The variable-step BDF3 formula 
$D_3v^n:=\bra{\Pi_{n,3}v}'(t_n)$ for $n\ge3$,
\begin{align}\label{eq: BDF3 formula}
D_3v^n=&\,d_0(r_n,r_{n-1})\partial_{\tau} v^{n}
+d_1(r_n,r_{n-1})\partial_{\tau} v^{n-1}+d_2(r_n,r_{n-1})\partial_{\tau} v^{n-2}.
\end{align}
As seen, the formula \eqref{eq: BDF3 formula} also represents
a variable-order variable-step BDF formula from the first-order to third-order.
The theoretical results on \eqref{eq: BDF3 formula} are naturally valid
for any combinations of the BDF-k formulas
if the time-step ratios are allowed to be zero in the associated coefficients.

We consider the variable-step BDF-$\rmk$ time-stepping scheme for the ODE problem,
\begin{align}\label{eq: BDF-k time-stepping}
	D_{\rmk}v^n=f(t_n,v^n)\quad\text{for $\rmk\le n\le N$.}
\end{align}
Without losing the generality, assume that the starting values $u^1$, $\cdots$, $u^{\rmk-1}$
are available by some other high-order approaches, such as Runge-Kutta methods.
For the stability analysis, assume further that the perturbed solution $\bar{v}^n$
solves the equation $D_{\rmk}\bar{v}^n=f(t_n,\bar{v}^n)+\varepsilon^n$
with a bounded sequence $\{\varepsilon^n\}$.
Then the difference $\tilde{v}^i=\bar{v}^i-v^i$ solves
\begin{align}\label{eq: perturbed equation}
	D_{\rmk}\tilde{v}^i=f(t_i,\bar{v}^i)-f(t_i,v^i)+\varepsilon^i\quad\text{for $\rmk\le i\le N$.}
\end{align}

These well-known stiff solvers have been tested to be practically valuable for differential-algebraic problems
\cite{DimarcoPareschi:2017,HairerNorsettWanner:1992,HairerWanner:2002}
and for the hyperbolic systems with multiscale relaxation
\cite{AlbiDimarcoPareschi:2020}. Also,
the variable-step versions are computationally efficient in capturing the multi-scale time behaviors
\cite{CalvoGrandeGrigorieff:1990,DeCariaGuzelLaytonLi:2021,HairerWanner:2002,LiaoJiZhang:2020pfc,LiaoJiWangZhang:2021,
	LiaoTangZhou:2020bdf2,WangRuuth:2008} via adaptive time-stepping strategies.
On the other hand, the stability and convergence analysis of variable-step BDF methods are difficult and
remain incomplete, cf. \cite{CalvoGrandeGrigorieff:1990,DeCariaGuzelLaytonLi:2021,Grigorieff:1983,
	GrigorieffPaes-Leme:1984,GuglielmiZennaro:2001,HairerNorsettWanner:1992,WangRuuth:2008},
because they involve multiple degrees of freedom (independent time-step sizes).

One of the main defects in the existing theory is the severe mesh condition required for stability.
For examples, the zero-stability of variable-step BDF2 method requires a step-ratio condition 
$0<r_k\le 1+\sqrt{2}\approx 2.414$, 
see \cite{Grigorieff:1983,GrigorieffPaes-Leme:1984}. 
By considering a specific time grid with constant step-ratio $r_k=r$, it was also shown in \cite{Grigorieff:1983} that
the zero-stability of variable-step BDF3 method requires the maximum step-ratio limit 
$R_3<\brat{1+\sqrt{5}}/{2}\approx1.618$. By using an elliptic type norm, 
Calvo, Grande and Grigorieff \cite{CalvoGrandeGrigorieff:1990}
gave the zero-stability constraint $r_k\le 1.476$ in 1990.
Twenty years ago, Guglielmi and Zennaro \cite{GuglielmiZennaro:2001} applied a spectral radius approach
to find the step-ratio condition $r_k\le 1.501$,
which seems to be the best result for the variable-step BDF3 method.

We aim to relax the existing mesh conditions for the perturbed stability of variable-step BDF2 and 
BDF3 methods according to the following definition.

\begin{definition}[\textbf{Stability}]\label{Def:stability}
	Assume that the nonlinear function $f(t,v)$ is a Lipschitz-continuous function 
	with the Lipschitz constant $L_f>0$. The $\rmk$-step time-stepping \eqref{eq: BDF-k time-stepping}
	with unequal time-steps	is stable if there exists a real number $\tau_0$ and a fixed constant $C$ (independent of the step sizes $\tau_n$) such that
	the numerical solution of \eqref{eq: perturbed equation} fulfills 
	\begin{align*}
		\max_{\rmk\le n\le N }\absb{\tilde{v}^{i}}\le C\bra{\sum_{i=0}^{\rmk-1}\absb{\tilde{v}^{i}}+\tau\sum_{i=1}^{\rmk-1}\absb{\partial_{\tau}\tilde{v}^{i}}
			+\max_{\rmk\le i\le N}\absb{\varepsilon^i}}\quad\text{for all $\tau_n\le \tau_0$.}
		\end{align*}
\end{definition}

This definition is practically reasonable in the sense that one always uses some high-order scheme
to start the $\rmk$-step time-stepping \eqref{eq: BDF-k time-stepping}.
In the uniform time-step case, it reduces 
into the classical zero-stability, cf. \cite[Definition 2.2]{CrouzeixLisbona:1984}, 
	\begin{align*}
	\max_{\rmk\le n\le N }\absb{\tilde{v}^{i}}\le C\bra{\sum_{i=0}^{\rmk-1}\absb{\tilde{v}^{i}}
		+\max_{\rmk\le i\le N}\absb{\varepsilon^i}}\quad\text{for all $\tau_n\le \tau_0$.}
\end{align*}
%On the non-uniform time grids, the zero stability method is stable, 
%while a stable method is not necessary zero-stable.
Definition \ref{Def:stability} would be as useful as 
the original zero-stability for the numerical analysis 
of variable-step, multi-step methods.
Actually, it was proven in \cite{LiaoZhang:2021} that the step-ratio condition $r_k< 3.561$
is sufficient to the stability of variable-step BDF2 method 
for linear parabolic problem. In a recent work \cite{LiaoJiWangZhang:2021}, 
this step-ratio condition was slightly improved to a new one, $r_k< 4.864$.
These step-ratio conditions  
are weaker than the classical zero-stability restriction $r_k<2.414$, but both of them are
are sufficient for the stability of BDF2 method
in the sense of Definition \ref{Def:stability}.

By developing a new framework with discrete orthogonal convolution (DOC) kernels,
we apply Definition \ref{Def:stability} to 
examine the ratio-step conditions for the stability of BDF2 and BDF3 methods in 
solving the ODE problems.
The main contribution is two-fold:
\begin{itemize}
  \item Theorem \ref{thm: BDF2 stability} shows that
  the BDF2 method is stable on arbitrary time grids.  
  \item Theorem \ref{thm: BDF3 stability} shows that the variable-step BDF3 method is stable
  if the adjacent step ratios $r_k< R_3,$ where $R_3\approx2.553$ is the unique positive root of
$$R_3^6+R_3^5-4 R_3^4-8 R_3^3-10 R_3^2-6R_3-2=0.$$
\end{itemize}

We start from a novel viewpoint by writing the BDF-$\rmk$ ($\rmk=2,3$) formula as
a convolution summation of local difference quotients
(this viewpoint is quite different from that in previous works
 \cite{LiaoJiWangZhang:2021,LiaoZhang:2021})
\begin{align}\label{eq: alter BDFk formula}
D_{\rmk}v^n:=\sum_{j=1}^nd^{(\rmk,n)}_{n-j}\partial_{\tau} v^j
\quad\text{for $n\ge\rmk$,}
\end{align}
where the first superscript index $\rmk$ in the discrete kernels $d^{(\rmk,n)}_{n-j}$
denotes the order of method. 
%We always need some starting procedures to initiate the BDF-$\rmk$ methods.
%Without losing generality,
%we assume the starting values are available for the BDF-$\rmk$ methods.
From the BDF2 formula \eqref{eq: BDF2 formula} with $n\ge2$, we have
\begin{align}
&d^{(2,n)}_j:=d_j(r_n,0)\;\;\text{for $j=0,1$,}\quad\text{and}\quad
d^{(2,n)}_j:=0\quad\text{for $n\ge j+1\ge 3$}.\label{eq: BDF2 kernels}
\end{align}
Similarly, from the BDF3 formula \eqref{eq: BDF3 formula} with $n\ge3$, we have
\begin{align}
&d^{(3,n)}_j:=d_j(r_n,r_{n-1})\;\;\text{for $j=0,1,2$,}\quad\text{and}\quad
d^{(3,n)}_j:=0\quad\text{for $n\ge j+1\ge 4$}.\label{eq: BDF3 kernels}
\end{align}

Here and hereafter, assume that the summation $\sum_{k=i}^{j}\cdot$ to be zero
and the product $\prod_{k=i}^{j}\cdot$ to be one if the index $i>j$.
As for the BDF-$\rmk$ kernels $d^{(\rmk,n)}_{n-j}$ with any fixed indexes $n$,
we recall a class of discrete orthogonal convolution
(DOC) kernels $\big\{\vartheta_{n-j}^{(\rmk,n)}\big\}_{j=\rmk}^n$ by a recursive procedure, also see \cite{LiaoZhang:2021,LiaoJiWangZhang:2021},
\begin{align}\label{eq: BDFk-DOC procedure}
\vartheta_{0}^{(\rmk,n)}:=\frac{1}{d^{(\rmk,n)}_{0}}\quad\text{and}\quad
\vartheta_{n-j}^{(\rmk,n)}:=-\frac{1}{d^{(\rmk,j)}_{0}}
\sum_{i=j+1}^{n}\vartheta_{n-i}^{(\rmk,n)}d^{(\rmk,i)}_{i-j}\quad\text{for $\rmk\leq j\le n-1$.}
\end{align}
Obviously, the DOC kernels $\vartheta_{n-j}^{(\rmk,n)}$ satisfy the following discrete orthogonality identity
\begin{align}\label{eq: BDFk-DOC identity}
\sum_{i=j}^{n}\vartheta_{n-i}^{(\rmk,n)}d^{(\rmk,i)}_{i-j}\equiv\delta_{nj}
\quad\text{for any $\rmk\leq j\le n$,}
\end{align}
where $\delta_{nj}$ is the Kronecker delta symbol with $\delta_{nj}=0$ if $j\neq n$.
Furthermore, with the
identity matrix $I_{m\times m}$ ($m:=n-\rmk+1$), the 
above discrete orthogonality identity \eqref{eq: BDFk-DOC identity} also implies
$$\mathbf{\Theta}_\rmk\mathbf{D}_{\rmk}=I_{m\times m},$$
where the two $m\times m$ matrices $\mathbf{D}_{\rmk}$ 
and $\mathbf{\Theta}_\rmk$ are defined by
\[
\mathbf{D}_{\rmk}:=
\left(
\begin{array}{cccc}
	d_{0}^{(\rmk,\rmk)}  &                  &  & \\
	d_{1}^{(\rmk,\rmk+1)}  &d_{0}^{(\rmk,\rmk+1)}  &  & \\
	\vdots           &\vdots           &\ddots  &\\
	d_{n-\rmk}^{(\rmk,n)}&d_{n-\rmk-1}^{(\rmk,n)}&\cdots  &d_{0}^{(\rmk,n)}  \\
\end{array}
\right)\;\;\text{and}\;\;
\mathbf{\Theta}_\rmk:=
\left(
\begin{array}{cccc}
	\vartheta_{0}^{(\rmk,\rmk)}  &                  &  & \\
	\vartheta_{1}^{(\rmk,\rmk+1)}  &\vartheta_{0}^{(\rmk,\rmk+1)}  &  & \\
	\vdots           &\vdots           &\ddots  &\\
	\vartheta_{n-\rmk}^{(\rmk,n)}&\vartheta_{n-\rmk-1}^{(\rmk,n)}&\cdots  &\vartheta_{0}^{(\rmk,n)}  \\
\end{array}
\right).
\]
Obviously, one has
$$\mathbf{D}_{\rmk}\mathbf{\Theta}_\rmk=I_{m\times m},$$
which implies the following mutual orthogonality identity
\begin{align}\label{eq: mutual BDFk-DOC identity}
	\sum_{i=j}^{n}d_{n-i}^{(\rmk,n)}\vartheta^{(\rmk,i)}_{i-j}\equiv\delta_{nj}
	\quad\text{for any $\rmk\leq j\le n$.}
\end{align}
We will use this identity \eqref{eq: mutual BDFk-DOC identity} 
to study the decaying property of DOC kernels $\vartheta_{n-j}^{(\rmk,n)}$.

We will derive an equivalent convolution form of the difference equation \eqref{eq: perturbed equation}. 
By exchanging the summation order and using \eqref{eq: BDFk-DOC identity}, one has
\begin{align}\label{eq: BDFk-DOC transform}
\sum_{i=\rmk}^{n}\vartheta_{n-i}^{(\rmk,n)}D_{\rmk}\tilde{v}^i
=&\,\sum_{i=\rmk}^{n}\vartheta_{n-i}^{(\rmk,n)}
\sum_{j=1}^{\rmk-1}d^{(\rmk,i)}_{i-j}\partial_{\tau}  v^j
+\sum_{i=\rmk}^{n}\vartheta_{n-i}^{(\rmk,n)}
\sum_{j=\rmk}^id^{(\rmk,i)}_{i-j}\partial_{\tau}  v^j\nonumber\\
=&\,\sum_{j=1}^{\rmk-1}\partial_{\tau}v^j
\sum_{i=\rmk}^n\vartheta_{n-i}^{(\rmk,n)}d^{(\rmk,i)}_{i-j}+
\sum_{j=\rmk}^{n}\partial_{\tau}v^j
\sum_{i=j}^n\vartheta_{n-i}^{(\rmk,n)}d^{(\rmk,i)}_{i-j}\nonumber\\
=&\,\mathcal{I}_{\rmk}^n[v]+\partial_{\tau}  v^n\qquad\text{for $n\ge\rmk$,}
\end{align}
where $\mathcal{I}_{\rmk}^n[v]$ represents the starting effect on the numerical solution at $t_n$,
\begin{align}\label{eq: initial effect BDFk-DOC}
	\mathcal{I}_{\rmk}^n[v]:=&\,\sum_{j=1}^{\rmk-1}\partial_{\tau}v^j
	\sum_{i=\rmk}^n\vartheta_{n-i}^{(\rmk,n)}d^{(\rmk,i)}_{i-j}\qquad\text{for $n\ge\rmk$.}
\end{align}
Multiplying both sides of the difference equation \eqref{eq: perturbed equation} by the DOC kernels
$\vartheta_{n-i}^{(\rmk,n)}$, and summing $i$ from $\rmk$ to $n$,
we apply \eqref{eq: BDFk-DOC transform} to get 
\begin{align}\label{eq: DOC perturbed equation}
	\partial_{\tau}\tilde{v}^n
	=-\mathcal{I}_{\rmk}^n[\tilde{v}]+\sum_{i=\rmk}^n\vartheta_{n-i}^{(\rmk,n)}
	\kbra{f(t_i,\bar{v}^i)-f(t_i,v^i)+\varepsilon^i}
	\quad\text{for $\rmk\le n\le N$.}
\end{align}
In next two sections,
we will present the stability analysis for the variable-step BDF2 and BDF3 methods, respectively,
via this discrete convolution form \eqref{eq: DOC perturbed equation}
of the perturbed equation. 
Numerical experiments on the graded and random time meshes are included 
in Section 4 to support our theory.
Some concluding remarks end this article.

\section{Unconditional stability of BDF2 method}
\setcounter{equation}{0}

%\newpage
For the BDF2 method, the associated DOC kernels are positive and decay exponentially.
\begin{lemma}\label{lem: BDF2 orthogonal formula}
The DOC kernels $\vartheta_{n-j}^{(2,n)}$ 
in \eqref{eq: BDFk-DOC procedure} have an explicit formula
\begin{align*}%%\label{eq: BDF2 orthogonal formula}
\vartheta_{n-j}^{(2,n)}=\frac{1}{d^{(2,j)}_{0}}
\prod_{i=j+1}^n\frac{r_{i}}{1+2r_{i}}>0\quad\text{ for $2\leq j\le n$},
\end{align*}
and satisfy %$\sum_{j=1}^{n}\vartheta_{n-j}^{(2,n)}=1$ for $n\ge2$.
\begin{align*}
\sum_{j=2}^{n}\vartheta_{n-j}^{(2,n)}=1-\prod_{i=2}^n\frac{r_{i}}{1+2r_{i}}<1\quad\text{for $n\ge2$.}
\end{align*}
\end{lemma}
\begin{proof}
For any fixed $j$, by taking $n=j$ in the identity
\eqref{eq: mutual BDFk-DOC identity}, one has
$\vartheta_{0}^{(2,j)}=1/d^{(2,j)}_{0}$
for $j\ge2$.
According to the definition \eqref{eq: BDF2 kernels}, one has
$d^{(2,n)}_{n-i}=0$ for $n\ge i+2$.
The identity \eqref{eq: mutual BDFk-DOC identity} gives
\begin{align*}
\vartheta_{m-j}^{(2,m)}=
-\frac{d^{(2,{m})}_{1}}{d^{(2,m)}_{0}}\vartheta_{m-1-j}^{(2,m-1)}
=-\frac{d_{1}(r_m,0)}{d_{0}(r_m,0)}\vartheta_{m-1-j}^{(2,m-1)}
=\frac{r_m}{1+2r_m}\vartheta_{m-1-j}^{(2,m-1)}
\end{align*}
for $m\ge j+1\ge2$. It yields immediately
\begin{align*}%%\label{eq: BDF2 orthogonal formula}
\vartheta_{n-j}^{(2,n)}
=&\,\vartheta_{n-2-j}^{(2,n-2)}\prod_{i=n-1}^n\frac{r_{i}}{1+2r_{i}}
=\cdots
=\vartheta_{0}^{(2,j)}\prod_{i=j+1}^n\frac{r_{i}}{1+2r_{i}}\quad\text{ for $2\leq j\le n$}.
\end{align*}
It leads to the first result.
%Moreover, the variable-step BDF2 method \eqref{eq: alter BDFk formula}
%of $\rmk=2$ is exact for the linear polynomial $v=t$.
By taking $v^n=t_n$ in \eqref{eq: alter BDFk formula}, 
one has $\partial_\tau t_n=1$ and $D_2t_n=1$ for $n\ge2$.
Thus applying the equality \eqref{eq: BDFk-DOC transform} with $\rmk=2$ yields
\begin{align}\label{lem: DOC summation BDF2}
\sum_{j=2}^{n}\vartheta_{n-j}^{(2,n)}
=&\,\sum_{j=2}^{n}\vartheta_{n-j}^{(2,n)}D_2t_j
=(\partial_{\tau}t_1)\sum_{i=2}^n\vartheta_{n-i}^{(2,n)}d^{(2,i)}_{i-1}+\partial_{\tau}t_n\nonumber\\
=&\,1+\vartheta_{n-2}^{(2,n)}d^{(2,2)}_{1}
=1-\prod_{i=2}^n\frac{r_{i}}{1+2r_{i}}<1\quad\text{for $n\ge1$,}
\end{align}
where the explicit formula of $\vartheta_{n-2}^{(2,n)}$ 
and the definition \eqref{eq: BDF2 kernels} were used.
It completes the proof.
\end{proof}

\begin{theorem}\label{thm: BDF2 stability}
If $\tau\le 1/(4L_f)$, the BDF2 solution of \eqref{eq: perturbed equation}  satisfies
\begin{align*}
\absb{\tilde{v}^{n}}\le 2\exp(4L_ft_{n-1})\braB{\absb{\tilde{v}^{1}}+2\tau\absb{\partial_{\tau}\tilde{v}^{1}}
+2t_{n}\max_{2\le i\le n}\absb{\varepsilon^i}}\quad\text{for $2\le n\le N$.}
\end{align*}
Thus the variable-step BDF2 scheme is stable on arbitrary time meshes.
\end{theorem}
\begin{proof} Take $\rmk=2$. 
Multiplying both sides of \eqref{eq: DOC perturbed equation} with $2\tau_n\tilde{v}^n$, one
applies Lemma \ref{lem: BDF2 orthogonal formula} (the positivity and boundedness of DOC kernels) 
to obtain that
\begin{align*}
\absb{\tilde{v}^n}^2-&\,\absb{\tilde{v}^{n-1}}^2+\absb{\tau_n\partial_{\tau}\tilde{v}^n}^2
=-2\tau_n\tilde{v}^n\mathcal{I}_{2}^n[\tilde{v}]
+2\tau_n\tilde{v}^n\sum_{i=2}^n\vartheta_{n-i}^{(2,n)}\kbra{f(t_i,\bar{v}^i)-f(t_i,v^i)+\varepsilon^i}\\
\le&\,2\tau_n\absb{\tilde{v}^n}\absb{\mathcal{I}_{2}^n[\tilde{v}]}
+2\tau_n\absb{\tilde{v}^n}\sum_{i=2}^n\vartheta_{n-i}^{(2,n)}\absb{f(t_i,\bar{v}^i)-f(t_i,v^i)}
+2\tau_n\absb{\tilde{v}^n}\max_{1\le i\le n}\absb{\varepsilon^i}\sum_{i=2}^n\vartheta_{n-i}^{(2,n)}\\
\le&\,2\tau_n\absb{\tilde{v}^n}\absb{\mathcal{I}_{2}^n[\tilde{v}]}
+2L_f\tau_n\absb{\tilde{v}^n}\sum_{i=2}^n\vartheta_{n-i}^{(2,n)}\absb{\tilde{v}^i}
+2\absb{\tilde{v}^n}\tau_n\max_{1\le i\le n}\absb{\varepsilon^i},
\end{align*}
where the Lipschitz continuous property of the nonlinear function $f$ has been used.
By dropping the nonnegative term at the left side and summing $n$ from $n=2$ to $m$, we have
\begin{align*}
\absb{\tilde{v}^m}^2\le&\,\absb{\tilde{v}^{1}}^2
+2\sum_{j=2}^m\tau_j\absb{\tilde{v}^j}\absb{\mathcal{I}_{2}^j[\tilde{v}]}
+2L_f\sum_{j=2}^m\tau_j\absb{\tilde{v}^j}\sum_{i=2}^j\vartheta_{j-i}^{(2,j)}\absb{\tilde{v}^i}
+2\sum_{j=2}^m\tau_j\absb{\tilde{v}^j}\max_{1\le i\le j}\absb{\varepsilon^i}.
\end{align*}
Let $m_0$ ($1\le m_0\le m$) be an integer such that $\absb{\tilde{v}^{m_0}}=\max_{1\le k\le m}\absb{\tilde{v}^{k}}$.
Now we take $m:=m_0$ in the above inequality and get
\begin{align*}
\absb{\tilde{v}^{m_0}}^2\le&\,\absb{\tilde{v}^{1}}\absb{\tilde{v}^{m_0}}
+2\absb{\tilde{v}^{m_0}}\sum_{j=2}^{m_0}\tau_j\absb{\mathcal{I}_{2}^j[\tilde{v}]}
+2L_f\absb{\tilde{v}^{m_0}}\sum_{j=2}^{m_0}\tau_j\absb{\tilde{v}^j}
+2\absb{\tilde{v}^{m_0}}\max_{2\le i\le m_0}\absb{\varepsilon^i}\sum_{j=2}^{m_0}\tau_j.
\end{align*}
It leads to
\begin{align}\label{thmproof: BDF2 stability}
\absb{\tilde{v}^{m}}\le \absb{\tilde{v}^{m_0}}\le\absb{\tilde{v}^{1}}
+2\sum_{j=2}^m\tau_j\absb{\mathcal{I}_{2}^j[\tilde{v}]}
+2L_f\sum_{j=2}^{m}\tau_j\absb{\tilde{v}^j}
+2t_{m}\max_{2\le i\le m}\absb{\varepsilon^i}\quad\text{for $2\le m\le N$.}
\end{align}
By using Lemma \ref{lem: BDF2 orthogonal formula}, the starting error $\mathcal{I}_{2}^n[\tilde{v}]$ 
in \eqref{eq: initial effect BDFk-DOC} reads
\begin{align*}
	\mathcal{I}_{2}^n[\tilde{v}]=\partial_{\tau}\tilde{v}^1\sum_{i=2}^n\vartheta_{n-i}^{(2,n)}d^{(2,i)}_{i-1}
	=\vartheta_{n-2}^{(2,n)}d^{(2,2)}_{1}(\partial_{\tau}\tilde{v}^1)
	=-(\partial_{\tau}\tilde{v}^1)\prod_{i=2}^n\frac{r_{i}}{1+2r_{i}}
	\quad\text{for $n\ge2$}
\end{align*}
such that
\begin{align*}
	\sum_{j=2}^m\tau_j\absb{\mathcal{I}_{2}^j[\tilde{v}]}
	\le\tau\absb{\partial_{\tau}\tilde{v}^1}\sum_{j=2}^m\frac{1}{2^{j-1}}
	\le\tau\absb{\partial_{\tau}\tilde{v}^1}
	\quad\text{for $m\ge2$.}
\end{align*}
It follows from \eqref{thmproof: BDF2 stability} that
\begin{align*}
	\absb{\tilde{v}^{n}}\le\absb{\tilde{v}^{1}}
	+2\tau\absb{\partial_{\tau}\tilde{v}^1}
	+2L_f\sum_{j=2}^{n}\tau_j\absb{\tilde{v}^j}
	+2t_{n}\max_{2\le i\le n}\absb{\varepsilon^i}\quad\text{for $2\le n\le N$.}
\end{align*}
Assuming that $\tau\le 1/(4L_f)$, one has
\begin{align*}
\absb{\tilde{v}^{n}}\le 2\absb{\tilde{v}^{1}}+4\tau\absb{\partial_{\tau}\tilde{v}^1}
+4L_f\sum_{j=2}^{n-1}\tau_j\absb{\tilde{v}^j}
+4t_{n}\max_{2\le i\le n}\absb{\varepsilon^i}\quad\text{for $2\le n\le N$.}
\end{align*}
The standard discrete Gr\"{o}nwall inequality, e.g. \cite[Lemma 3.1]{LiaoZhang:2021}, completes the proof.
\end{proof}

%\begin{remark}
%	The unconditional zero-stability 
%	can be understanded
%	in a straightforward manner by considering the homogeneous equation
%	$D_{2}v^n=0$ or
%	\begin{align*}
%		\partial_{\tau} v^n=-\frac{d^{(2,n)}_{1}}{d^{(2,n)}_{0}}\partial_{\tau} v^{n-1}
%		=\frac{r_n}{1+2r_n}\partial_{\tau} v^{n-1}
%		\quad\text{for $n\ge2$.}
%	\end{align*}
%	It leads to
%	\begin{align*}
%		\partial_{\tau} v^n=\partial_{\tau} v^{n-2}\prod_{i=n-1}^n\frac{r_i}{1+2r_i}=\cdots
%		=\partial_{\tau} v^{1}\prod_{i=2}^n\frac{r_i}{1+2r_i}\quad\text{for $n\ge2$.}
%	\end{align*}
%	It follows that
%	\begin{align*}
%		v^n=v^{n-1}+\tau_n(\partial_{\tau} v^{1})\prod_{i=2}^n\frac{r_i}{1+2r_i}\quad\text{for $n\ge2$}
%	\end{align*}
%	and then
%	\begin{align*}
%		v^n=v^{1}+(\partial_{\tau} v^{1})\sum_{j=2}^n\tau_j\prod_{i=2}^j\frac{r_i}{1+2r_i}\quad\text{for $n\ge2$.}
%	\end{align*}
%	The unconditional zero-stability is checked again 
%	since $|v^n|\le |v^1|+\tau|\partial_{\tau} v^{1}|$ for $n\ge2$.	
%    On the other hand, our result in Theorem \ref{thm: BDF2 stability} is proved by the discrete energy technique,
%    which is extendable to other BDF methods. 	
%\end{remark}

\section{Stability analysis of BDF3 method}
\setcounter{equation}{0}

\subsection{Decaying of DOC kernels}

Note that, the variable-step BDF3 method \eqref{eq: alter BDFk formula}
of $\rmk=3$ is exact for the linear polynomial $v=t$.
Taking $v^n=t_n$ in \eqref{eq: alter BDFk formula},
one can find that $D_3t_n=1$ for $n\ge3$. As done in \eqref{lem: DOC summation BDF2},
we can apply the discrete equality \eqref{eq: BDFk-DOC transform} 
with \eqref{eq: initial effect BDFk-DOC} to get
\begin{align*}
	\sum_{j=3}^{n}\vartheta_{n-j}^{(3,n)}
	=&\,\sum_{j=1}^{2}\sum_{i=3}^n\vartheta_{n-i}^{(3,n)}d^{(3,i)}_{i-j}+1\quad\text{for $n\ge3$.}
\end{align*}
However, it does not provide enough information for the subsequent stability analysis 
because no explicit formulas of DOC kernels $\vartheta_{n-j}^{(3,n)}$ are available.
Furthermore, the DOC kernels are not always positive,
see Remark \ref{remark: DOC values-BDF3} below, we turn to bound
$\sum_{j=3}^{n}\absb{\vartheta_{n-j}^{(3,n)}}$ for any time-levels $t_n$
by imposing certain step-ratio constraint. 

To do so,
introduce the modified DOC kernels
\begin{align}\label{lem: auxiliary DOC BDF3}
\widehat{\vartheta}_{i-j}^{(3,i)}:=\vartheta_{i-j}^{(3,i)}d^{(3,j)}_{0}\quad\text{for $i\ge j\ge 3$.}
\end{align}
The discrete identity \eqref{eq: mutual BDFk-DOC identity} gives
\begin{align}\label{eq: auxiliary BDFk-DOC identity}
\sum_{i=j}^{n}d_{n-i}^{(3,n)}\widehat{\vartheta}^{(3,i)}_{i-j}\equiv\delta_{nj}d^{(3,j)}_{0}
\quad\text{for $n\geq j\ge 3$.}
\end{align}
For any fixed index $j\ge3$, by taking $n=j$
and $n=j+1$ in the identity \eqref{eq: auxiliary BDFk-DOC identity}, respectively, one can derive that
\begin{align}\label{eq: DOC-BDF3 initial data}
\widehat{\vartheta}^{(3,j)}_{0}=1
\quad\text{and}\quad
\widehat{\vartheta}^{(3,j+1)}_{1}=-\frac{d_{1}^{(3,j+1)}}{d_{0}^{(3,j+1)}}
\quad\text{for $j\ge3$.}
\end{align}
According to the definition \eqref{eq: BDF3 kernels}, one has
$d^{(3,n)}_{n-i}=0$ for $n\ge i+3$.
The identity \eqref{eq: mutual BDFk-DOC identity} gives
\begin{align*}
\widehat{\vartheta}_{m-j}^{(3,m)}+\frac{d^{(3,{m})}_{1}}{d^{(3,m)}_{0}}\widehat{\vartheta}_{m-1-j}^{(3,m-1)}
+\frac{d^{(3,{m})}_{2}}{d^{(3,m)}_{0}}\widehat{\vartheta}_{m-2-j}^{(3,m-2)}=0
\quad\text{for $m\ge j+2\ge5$,}
\end{align*}
or the difference equation
\begin{align}\label{eq: DOC-BDF3 difference equation}
\widehat{\vartheta}_{m-j}^{(3,m)}-\alpha_m\widehat{\vartheta}_{m-1-j}^{(3,m-1)}
+\beta_m\widehat{\vartheta}_{m-2-j}^{(3,m-2)}=0
\quad\text{for $m\ge j+2\ge5$,}
\end{align}
where, by using \eqref{funs: d0}-\eqref{funs: d2} and
\eqref{eq: alpha auxiliary function}-\eqref{eq: beta auxiliary function},
the coefficients $\alpha_m$ and $\beta_m$ are defined by
\begin{align}\label{eq: DOC-BDF3 coefficients}
\alpha_m:=-\frac{d^{(3,{m})}_{1}}{d^{(3,m)}_{0}}=\alpha(r_m,r_{m-1})
\quad\text{and}\quad
\beta_m:=\frac{d^{(3,{m})}_{2}}{d^{(3,m)}_{0}}=\beta(r_m,r_{m-1})
\quad\text{for $m\ge 3$.}
\end{align}

To evaluate the asymptotic behaviors of the DOC kernels $\widehat{\vartheta}_{m-j}^{(3,m)}$,
introduce the following  vector and companion matrix
\begin{align}\label{eq: DOC-BDF3 matrix}
\vec{u}_m=\left(
            \begin{array}{c}
              \widehat{\vartheta}_{m-j}^{(3,m)} \\[2ex]
              \widehat{\vartheta}_{m-1-j}^{(3,m-1)} \\
            \end{array}
          \right)
\quad\text{and}\quad
A_m:=\left(
      \begin{array}{cc}
        \alpha_m & -\beta_m \\[2ex]
        1 & 0 \\
      \end{array}
    \right).
\end{align}
For any fixed index $j\ge1$, the difference equation \eqref{eq: DOC-BDF3 difference equation} reads
\begin{align*}
\vec{u}_m=A_m\vec{u}_{m-1}
\quad\text{for $m\ge j+2\ge5$,}
\end{align*}
or
\begin{align}\label{eq: DOC-BDF3 matrix equation}
\vec{u}_n=A_n\vec{u}_{n-1}=A_nA_{n-1}\vec{u}_{n-2}
=\cdots=\prod_{i=j+2}^nA_i\vec{u}_{j+1}
\quad\text{for $n\ge j+2\ge5$.}
\end{align}
The associated initial vector
$\vec{u}_{j+1}=\brat{\alpha_{j+1}, 1}^T$
according to \eqref{eq: DOC-BDF3 initial data}.

If all of step-ratios $0<r_k<\hat{R}_3\approx3.4405$ ($k\ge2$),
Lemmas \ref{lemma: alpha beta 1}-\ref{lemma: beta less than 1} show that
\begin{align}\label{eq: DOC-BDF3 coefficients property}
\beta_i<\alpha_i<1+\beta_i<2 \quad\text{for $i\ge 3$.}
\end{align}
 They imply that the roots of $\lambda^2-\alpha_i\lambda+\beta_i=0$ satisfy $\abs{\lambda}<1$.
 Thus the eigenvalues of $A_i$ have the module less than 1 and
 the spectral radius $\rho(A_i)<1$ for any $i\ge3$.
 In general, it is not sufficient to ensure the global decaying
 of DOC kernels $\widehat{\vartheta}_{n-j}^{(3,n)}$.

 It needs to examine
 the joint companion matrix $\prod_{i=j+2}^nA_i$ for $n\ge j+2\ge5$.
Let $\mynorm{\cdot}_{\infty}$ be the maximum
norm on the space $\mathbb{C}^2$ and let the same symbol $\mynorm{\cdot}_{\infty}$ denote also the
corresponding induced matrix norm.
We adopt the technique of Calvo \emph{et al} \cite{CalvoGrandeGrigorieff:1990}
by considering a complex constant $\mu$ with $\mathrm{Im}(\mu)\neq0$
and the following elliptic type norm
\begin{align}\label{eq: elliptic type norm}
\mynorm{A_i}_{H,\infty}:=\mynorm{H^{-1}A_iH}_{\infty}\quad\text{with}\quad
H:=\left(
      \begin{array}{cc}
        \mu & \bar{\mu} \\[2ex]
        1 & 1 \\
      \end{array}
    \right).
    \end{align}
Then one can bound the DOC kernels via
\begin{align}\label{eq: DOC-BDF3 norm}
\mynormb{\vec{u}_n}_{\infty}\leq
\mynormb{H}_{\infty}\prod_{i=j+2}^n\mynormb{A_i}_{H,\infty}\mynormb{H^{-1}}_{\infty}\mynormb{\vec{u}_{j+1}}_{\infty}
\quad\text{for $n\ge j+2\ge5$.}
\end{align}

To process the analysis, it is to determine a fixed parameter $\mu$.
We will choose a proper constant $\tilde{R}_3<\hat{R}_3\approx3.4405$ to ensure that
\begin{align}\label{eq: DOC-elliptic norm less than 1}
\mynormb{A_i}_{H,\infty}<1\quad\text{for $i\ge 3$,}
\end{align}
if all of step-ratios $r_k<\tilde{R}_{3}$ ($k\ge2$).
It is easy to derive that
\begin{align}\label{def: DOC-elliptic norm}
\mynormb{A_i}_{H,\infty}=\frac{\abs{\mu^2-\alpha_i\mu+\beta_i}
+\absb{\abs{\mu}^2-\alpha_i\mu+\beta_i}}{\abs{\mu-\bar{\mu}}}\quad\text{for $i\ge 3$.}
\end{align}
%where $p_i(\mu):=\mu^2-\alpha_i\mu+\beta_i$. 
Noticing that the following fact
\begin{align*}
	\absb{\abs{\mu}^2-\alpha_i\mu+\beta_i}^2-\absb{\mu^2-\alpha_i\mu+\beta_i}^2
	=\beta_i\absb{\mu-\bar{\mu}}^2,
\end{align*}
one  has
\begin{align}\label{eq: DOC-elliptic norm}
	\mynormb{A_i}_{H,\infty}=\frac{\abs{\mu^2-\alpha_i\mu+\beta_i}}{\abs{\mu-\bar{\mu}}}
	+\sqrt{\beta_i+\frac{\abs{\mu^2-\alpha_i\mu+\beta_i}^2}{\abs{\mu-\bar{\mu}}^2}}\quad\text{for $i\ge 3$.}
\end{align}
Then the requirement \eqref{eq: DOC-elliptic norm less than 1} is equivalent to
\begin{align}\label{eq: DOC-elliptic norm derive 2}
\frac{\absb{\mu^2-\alpha_i\mu+\beta_i}}{\abs{\mu-\bar{\mu}}}<\frac12(1-\beta_i)
\quad\text{for $i\ge 3$}.
\end{align}
Then, taking  into account the coefficient relationship \eqref{eq: DOC-BDF3 coefficients property}, we can find that
the condition \eqref{eq: DOC-elliptic norm derive 2} can be written equivalently in the following form,
also see \cite[Eq. (14)]{CalvoGrandeGrigorieff:1990},
\begin{align}\label{eq: DOC-equivalence condition}
\abs{\mu-\frac12\alpha_i-\frac{\imath}2 \sqrt{(1+\beta_i)^2-\alpha_i^2}}<\frac12(1-\beta_i)
\quad\text{for $i\ge 3$.}
\end{align}
Note that, the other branch requiring 
$\absb{\mu-\frac12\alpha_i+\frac{\imath}2\sqrt{(1+\beta_i)^2-\alpha_i^2}}<\frac12(1-\beta_i)$
is omitted here since it is equivalent to \eqref{eq: DOC-equivalence condition} 
by replacing the undetermined parameter $\mu$ with $\bar{\mu}$.
Since the coefficients $\alpha_i,\beta_i$ are dependent on the ratios $r_i$ and $r_{i-1}$,
The inequalities in \eqref{eq: DOC-equivalence condition}  define
a family of complex disks $\mathfrak{D}(r_i,r_{i-1})$
 centered at $\brab{\frac12\alpha_i,\frac12\sqrt{(1+\beta_i)^2-\alpha_i^2}}$
 with the radius $\frac12(1-\beta_i)$.

A heuristic analysis is considered firstly to obtain a rough estimate for the complex parameter $\mu$,
 while the mathematical proof is left to next subsection.
 To ensure \eqref{eq: DOC-equivalence condition}, one should choose a fixed $\mu_0$
 such that the intersection set
 $$\bigcap_{r_i,r_{i-1}\in[0, \tilde{R}_3]}\mathfrak{D}(r_i,r_{i-1})\quad\text{ is not empty.}$$
Reminding the increasing functions
\eqref{eq: alpha auxiliary function}-\eqref{eq: beta auxiliary function},
we see that the disk $\mathfrak{D}(0,0)$ centered at $\brab{0,\frac12}$
has the maximum radius $\frac12$, while the disk $\mathfrak{D}(\tilde{R}_3,\tilde{R}_3)$ centered at
$$\bra{\frac12\alpha(\tilde{R}_3,\tilde{R}_3),\frac12\sqrt{\kbra{1+\beta(\tilde{R}_3,\tilde{R}_3)}^2-\alpha^2(\tilde{R}_3,\tilde{R}_3)}}$$
has the minimum radius $\frac12\brab{1-\beta(\tilde{R}_3,\tilde{R}_3)}$. Obviously,
the largest value of $\tilde{R}_3$  may be determined from the fact that
the smallest disk $\mathfrak{D}(\tilde{R}_3,\tilde{R}_3)$ is tangential to the largest one $\mathfrak{D}(0,0)$
at the tangential point $\mu_0$.
This necessary condition holds if and only if
\begin{align*}
\frac14\alpha^2(\tilde{R}_3,\tilde{R}_3)
+\bra{\frac12\sqrt{\kbra{1+\beta(\tilde{R}_3,\tilde{R}_3)}^2-\alpha^2(\tilde{R}_3,\tilde{R}_3)}-\frac12}^2
=\frac14\kbrab{2-\beta(\tilde{R}_3,\tilde{R}_3)}^2\,.
\end{align*}
which leads to
\begin{align}\label{eq: guess step-ratio condition}
\alpha^2(\tilde{R}_3,\tilde{R}_3)+8\beta^2(\tilde{R}_3,\tilde{R}_3)-8\beta(\tilde{R}_3,\tilde{R}_3)=0,
\end{align}
or equivalently,
\begin{align*}
9 \tilde{R}_3^6-2 \tilde{R}_3^5-35 \tilde{R}_3^4-42 \tilde{R}_3^3-22 \tilde{R}_3^2-4 \tilde{R}_3+1=0.
\end{align*}
This polynomial equation has two positive roots $\tilde{R}_3\approx2.5808$
and $\tilde{R}_3\approx0.1304$
(it is checked that this situation also occurs in \cite{CalvoGrandeGrigorieff:1990}).
We throw away the small root $\tilde{R}_3\approx0.1304$ and choose the large one $\tilde{R}_3\approx2.5808$
with the corresponding tangential point $\tilde{\mu}_*\approx0.4979+0.5454\imath$.

\begin{remark}\label{remark: spectral radius approach}
It is to mention that, by following the joint spectral radius approach in \cite{GuglielmiZennaro:2001},
one can obtain a slightly improved step-ratio constraint $0<r_k<2.705$ to
ensure the boundedness of DOC kernels; however,
it seems difficult to obtain the desired decaying behavior of DOC kernels theoretically
from the extremal norm estimate of the family of companion matrices $\{A_i\,|\,i\ge3\}$.
\end{remark}

\begin{remark}\label{remark: DOC-equivalence condition}
The condition \eqref{eq: DOC-equivalence condition} is sufficient
for \eqref{eq: DOC-elliptic norm less than 1}, or $\mynorm{A_i}_{H,\infty}<1$;
but different choices of the parameter $\mu$ leads to different values of $\mynorm{A_i}_{H,\infty}$.
Consider the simple case with $r_{i-1}\equiv0$ for any $i\ge3$, corresponding to
the variable-step BDF2 scheme. By using the definition \eqref{eq: DOC-BDF3 coefficients}
together with \eqref{funs: d0}-\eqref{funs: d2}, one has
$\alpha_i=\frac{r_i}{1+2r_i}$, $\beta_i=0$ and the reduced class of complex discs
\begin{align*}
\mathfrak{D}(r_i,0)=\left\{\mu\,\left|\,\absB{\mu-\frac{1}{2}\alpha(r_i,0)-\frac{\imath}{2}
\sqrt{1-\alpha^2(r_i,0)}}<\frac12\right.\right\}
\quad\text{for $i\ge 3$.}
\end{align*}
Recalling the fact $(\alpha_i+\frac12)(\alpha_i-1)<0$ for any $r_i\ge0$, one can check that
 $$\abs{\frac{1}{4}-\frac{1}{2}\alpha(r_i,0)
 +\frac{\imath}{2}\braB{\frac{\sqrt{3}}{2}-\sqrt{1-\alpha^2(r_i,0)}}}<\frac{1}{2}.$$
That is, the distance between the center $\brab{\frac1{4},\frac{\sqrt{3}}4}$
of $\mathfrak{D}(\infty,0)$ and the center $\brab{\frac12\alpha_i,\frac12\sqrt{1-\alpha_i^2}}$
of $\mathfrak{D}(r_i,0)$ is always less than $\frac1{2}$.
Thus these complex discs are always have a common region for any $r_i>0$, that is,
 $\cap_{r_i>0}\mathfrak{D}(r_i,0)$ is not empty. There exists a fixed parameter $\mu$ such that
$\mynormb{A_i}_{H,\infty}<1$ for any $r_i>0$.
%By taking $\mu_0=\frac1{4}+\frac{\sqrt{3}}4\imath$ in \eqref{eq: DOC-elliptic norm}, 
%the corresponding elliptic norm takes the value of
%\begin{align*}
%	\mynormb{A_i}_{H,\infty}
%	=\frac{2\abs{\mu_0}\abs{\mu_0-\alpha_i}}{\abs{\mu_0-\bar{\mu}_0}}
%	=\frac{2}{\sqrt{3}}\sqrt{\brab{\frac{1}{4}-\alpha_i}^2+\frac{3}{16}}
%	\quad\text{for $r_i>0$.}
%\end{align*}
%Since $\alpha_i\in(0,\frac{1}{2})$ for any $r_i>0$, one gets $\mynormb{A_i}_{H,\infty}\in\brab{\frac12,\frac1{\sqrt{3}}}$.
By taking $\mu_{1}=\imath/2$ in \eqref{eq: DOC-elliptic norm}, one has
\begin{align*}
	\mynormb{A_i}_{H,\infty}
	=\frac{2\abs{\mu_{1}}\abs{\mu_{1}-\alpha_i}}{\abs{\mu_{1}-\bar{\mu}_{1}}}
	=\sqrt{\alpha_i^2+\frac{1}{4}}
	\quad\text{for $r_i>0$.}
\end{align*} 
Since $\alpha_i\in(0,\frac{1}{2})$ for any $r_i>0$, one gets 
$\mynorm{A_i}_{H,\infty}\in\brab{\tfrac{1}{2},\tfrac{1}{\sqrt{2}}}$.
On the other hand, taking $\mu_{2}=1/4+\imath/4$ in \eqref{eq: DOC-elliptic norm} arrives at
\begin{align*}
	\mynormb{A_i}_{H,\infty}
	=\frac{2\abs{\mu_{2}}\abs{\mu_{2}-\alpha_i}}{\abs{\mu_{2}-\bar{\mu}_{2}}}
	=\sqrt{2\brab{\frac{1}{4}-\alpha_i}^2+\frac{1}{8}}
	\quad\text{for $r_i>0$.}
\end{align*}
For any $r_i>0$, we have $\mynorm{A_i}_{H,\infty}\in\brab{\frac{\sqrt{2}}{4},\frac{1}{2}}$, which are less than 1/2.
Always, the estimate \eqref{eq: DOC-BDF3 norm} implies
the globally asymptotic decaying of DOC kernels, as
shown by Lemma \ref{lem: BDF2 orthogonal formula}.
\end{remark}

\subsection{Stability analysis}

To avoid the undesired equation \eqref{eq: guess step-ratio condition} which has two different positive roots,
and to simplify the subsequent mathematical derivations,
we fix the complex parameter
$$\mu=\mu_{*}:=1/2+\imath/2,$$
which is very close to the optimal tangential point $\tilde{\mu}_*\approx0.4979+0.5454\imath$ in the heuristic analysis.
Reminding the condition \eqref{eq: DOC-equivalence condition}, we determine the maximum step-ratio $R_3$ by
\begin{align*}
\abs{1-\alpha(R_3,R_3)+\imath \braB{1-\sqrt{(1+\beta(R_3,R_3))^2-\alpha(R_3,R_3)^2}}}=1-\beta(R_3,R_3),
\end{align*}
which leads to the polynomial equation
%$R_3^6+R_3^5-4 R_3^4-8 R_3^3-10 R_3^2-6 R_3-2=0$.
\begin{align*}
R_3^6+R_3^5-4 R_3^4-8 R_3^3-10 R_3^2-6 R_3-2=0.
\end{align*}
It has a unique positive root $R_3\approx2.553$. Now we prove the following result.

%\begin{align*}
%	\mynormb{A_i}_{H,\infty}=\frac{\abs{\mu^2-\alpha_i\mu+\beta_i}
%		+\absb{\abs{\mu}^2-\alpha_i\mu+\beta_i}}{\abs{\mu-\bar{\mu}}}\quad\text{for $i\ge 3$.}
%\end{align*}
%\begin{align*}
%	\mynormb{A_i}_{H,\infty}=\frac{\abs{\mu^2-\alpha_i\mu+\beta_i}
%		+\absb{\abs{\mu}^2-\alpha_i\mu+\beta_i}}{\abs{\mu-\bar{\mu}}}\quad\text{for $i\ge 3$.}
%\end{align*}

\begin{lemma}\label{lem: BDF3 orthogonal formula}
If the step ratios $r_k< R_3$ $(k\ge2)$, 
there exists a constant $C_3>0$,
independent of the time $t_n$,
such that the DOC kernels $\vartheta_{n-j}^{(3,n)}$ in \eqref{eq: BDFk-DOC procedure}
satisfy
\begin{align*}
\sum_{j=3}^{n}\absb{\vartheta_{n-j}^{(3,n)}}\le C_3
\quad\text{and}\quad
\sum_{j=i}^{n}\absb{\vartheta_{j-i}^{(3,j)}}\le C_3\quad\text{for $n\ge 3\;(i\ge3)$.}
\end{align*}
\end{lemma}
\begin{proof}
The definitions  \eqref{funs: d0} and \eqref{eq: BDF3 kernels}
give $d^{(3,j)}_{0}=d_0(r_j,r_{j-1})\ge1$ for $j\ge3$.
By using  \eqref{lem: auxiliary DOC BDF3}, \eqref{eq: DOC-BDF3 initial data} and
the coefficient bound \eqref{eq: DOC-BDF3 coefficients property}, one gets
$$\vartheta_{0}^{(3,j)}=\frac1{d^{(3,j)}_{0}}\le1\quad\text{and}\quad \vartheta_{1}^{(3,j+1)}
=\frac{\alpha_{j+1}}{d^{(3,j)}_{0}}\le \frac{2}{d^{(3,j)}_{0}}\le 2,
\quad\text{for $j\ge3$.}$$
That is, the first two kernels are always positive and bounded so that
\begin{align*}
\sum_{j=3}^n\absb{\vartheta_{n-j}^{(3,n)}}\le 3
\quad\text{and}\quad
\sum_{j=i}^n\absb{\vartheta_{j-i}^{(3,j)}}\le 3
\quad\text{for $n=3,4\; (i\ge3)$.}
\end{align*}

Now consider the general case $n\ge5$.
%$$\mynormb{H}_{\infty}=2, \qquad\mynormb{H^{-1}}_{\infty}=1+\frac{\sqrt{2}}2.$$
Taking $\mu=\mu^{*}$ in the condition \eqref{eq: DOC-elliptic norm derive 2}  yields
\begin{align*}
	\absb{2\beta_i-\alpha_i+\imath(1-\alpha_i)}<1-\beta_i	\quad\text{for $i\ge 3$,}
\end{align*}
or equivalently, 
\begin{align*}
	2\alpha_i^2+3\beta_i^2-4\alpha_i\beta_i-2\alpha_i+2\beta_i<0
	\quad\text{for $i\ge 3$.}
\end{align*}
Thanks to Lemma \ref{lemma: elliptic function}, they are valid if 
all of step-ratios $0< r_k<R_3$ $(k\ge2)$.
According to \eqref{eq: DOC-BDF3 matrix} and \eqref{eq: DOC-elliptic norm}, 
the elliptic norm $\mynormb{A_i}_{H,\infty}$ 
is a continuous function with respect to $r_i,r_{i-1}\in[\epsilon,R_3-\epsilon]$ for aribitrary small $\epsilon>0$.
Thus there exists a constant $\delta\in(0,1)$ such that
\begin{align*}
	\mynormb{A_i}_{H,\infty}\le\delta<1\quad\text{for $i\ge 3$.}
\end{align*}
By taking the parameter $\mu=\mu^{*}$
in \eqref{eq: elliptic type norm} and \eqref{eq: DOC-elliptic norm}, one has
$\mynormb{H}_{\infty}=2$ and $\mynormb{H^{-1}}_{\infty}=1+\frac{\sqrt{2}}2.$
Moreover, by \eqref{eq: alpha x-increasing}-\eqref{eq: alpha y-increasing},
one has $\mynormb{\vec{u}_{j+1}}_{\infty}\le \max\{\alpha(R_3,R_3),1\}=\alpha(R_3,R_3)$.
Thus, the estimate \eqref{eq: DOC-BDF3 norm} gives
\begin{align*}
	\absb{\vartheta_{n-j}^{(3,n)}}d^{(3,j)}_{0}\le\mynormb{\vec{u}_n}_{\infty}\leq
	c_R\delta^{n-j-1}
	\quad\text{for $n\ge j+2\ge5$,}
\end{align*}
where $c_R=(2+\sqrt{2})\alpha(R_3,R_3)$, and then
\begin{align}\label{eq: DOC-BDF3 decaying estimate}
	\absb{\vartheta_{n-j}^{(3,n)}}\leq c_R\delta^{n-j-1}
	\quad\text{for $n\ge j+2\ge5$.}
\end{align}
Therefore it follows that
\begin{align*}
	&\sum_{j=3}^{n}\absb{\vartheta_{n-j}^{(3,n)}}
	\le c_R\sum_{j=3}^{n}\delta^{n-j-1}<\frac{c_R}{\delta(1-\delta)}\quad\text{for $n\ge 5$,}\\
	&\sum_{j=i}^{n}\absb{\vartheta_{j-i}^{(3,j)}}
	\le c_R\sum_{j=i}^{n}\delta^{j-i-1}<\frac{c_R}{\delta(1-\delta)}\quad\text{for $n\ge i\ge3$.}
\end{align*}
We obtain the desired estimates by taking a constant $C_3:=\frac{c_R}{\delta(1-\delta)},$
which is independent of the time-level index $n$. This proof is complete.
\end{proof}

%\lan{\begin{remark}\label{remark: BDF3 orthogonal formula}
%		On one hand, we have a smaller value of $\delta_R$ 
%		if we set some smaller limit of step-ratios. 
%		Consider the step-ratios $r_k\le 2$, 
%		one can choose 
%		\begin{align*}
%			\delta_0=&\,\sqrt{\mathfrak{J}(\mu_0;1,2)}
%			+\sqrt{\mathfrak{J}(\mu_0;1,2)+\beta(1,2)}\approx0.6789,\\
%			\delta_*=&\,\sqrt{\mathfrak{J}(\mu_*;2,2)}
%			+\sqrt{\mathfrak{J}(\mu_*;2,2)+\beta(2,2)}\approx0.7540,
%		\end{align*} 
%		and then
%		$$c_R=2\sqrt{3}\max\{\alpha(2,2),1\}\approx3.6852,
%		\quad\delta_R=\max\{\delta_0,\delta_*\}<0.76.$$
%		By assuming $r_k\le 3/2$ $(k\ge2)$, 
%		the estimates \eqref{ieq: delta infty-elliptic norm 1},
%		\eqref{ieq: delta star-elliptic norm 1} and \eqref{eq: DOC-BDF3 decaying estimate} are also valid with 
%		\begin{align*}
%			\delta_0=&\,\sqrt{\mathfrak{J}(\mu_0;1,3/2)}
%			+\sqrt{\mathfrak{J}(\mu_0;1,3/2)+\beta(1,3/2)}\approx0.6312,\\
%			\delta_*=&\,\sqrt{\mathfrak{J}(\mu_*;1,0)}
%			+\sqrt{\mathfrak{J}(\mu_*;1,0)+\beta(1,0)}\approx0.7454,
%		\end{align*} 
%		and then
%		$$c_R=2\sqrt{3}\max\{\alpha(3/2,3/2),1\}=2\sqrt{3},\quad\delta_R=\max\{\delta_0,\delta_*\}<0.75.$$
%		On the other hand, the constant $C_3$ remains bounded 
%		as long as most of step-ratios $r_k<R_3$. Actually, if the element number of $\mathfrak{F}:=\{i_0\,|\,\mynorm{A_{i_0}}_{H,\infty}\ge1\}$ 
%		is finite or $\abs{\mathfrak{F}}\ll N$, 
%		the claimed results of Lemma \ref{lem: BDF3 orthogonal formula} are also valid. 
%		In other words, the imposed step-ratio condition $r_k<R_3$ can be relaxed to some extent.			
%		\end{remark}}

\begin{remark}\label{remark: elliptic type norm}
	The elliptic type norm in \eqref{eq: elliptic type norm}
	admits some other nonsingular matrices $H$.
	Here consider a simple case with
	$$H_0:=\Big(
	\begin{array}{cc}
		1 & 0 \\
		1 & 1 \\
	\end{array}
	\Big)\quad\text{such that}\quad
	H_0^{-1}A_iH_0=\left(
	\begin{array}{cc}
		\alpha_i-\beta_i & -\beta_i \\
		1-\alpha_i+\beta_i & \beta_i \\
	\end{array}
	\right)$$
	Let $R_{3,0}:=\frac{1}{3}+\frac{1}{3} \sqrt[3]{19+3 \sqrt{33}}+\frac{1}{3}\sqrt[3]{19-3 \sqrt{33}}\approx1.839$
	be the unique positive root of the cubic equation $R_{3,0}^3-R_{3,0}^2-R_{3,0}-1=0$.
	If the step-ratios satisfy $0<r_k<R_{3,0}$
	(it is also superior to the existing mesh conditions in the literature \cite{CalvoGrandeGrigorieff:1990,Grigorieff:1983,GuglielmiZennaro:2001}),
	one can follow the proofs of Lemmas \ref{lemma: alpha beta 1}-\ref{lemma: beta less than 1}
	to check that
	$$2\beta_i<\alpha_i<1\quad\text{ for $i\ge 3$.}$$
	It follows  that
	$$\mynormb{A_i}_{H_0,\infty}=\max\left\{\alpha_i,1-\alpha_i+2\beta_i\right\}<1
	\quad\text{for $i\ge 3$,}$$
	and then the DOC kernels $\vartheta_{n-j}^{(3,n)}$
	are globally decaying.
	Although the proof of Lemma \ref{lem: BDF3 orthogonal formula}
	is technically complex, the matrix $H$ in \eqref{eq: elliptic type norm} would be
	optimal in the sense that the companion matrix $A_i$ always has a pair of conjugate complex eigenvalues.
	Actually, the inequality $\alpha_i^2<4\beta_i$
	holds when the adjacent step ratios $r_i,r_{i-1}$ are close to 1,
	cf. Remark \ref{remark: DOC values-BDF3}.
\end{remark}

\begin{remark}\label{remark: DOC values-BDF3}
	Consider the uniform time-step $\tau_k=\tau$,
	the BDF3 kernels \eqref{eq: BDF3 kernels} give
	$$d_0^{(3,n)}=\frac{11}{6},\quad
	d_1^{(3,n)}=-\frac{7}{6},\quad
	d_2^{(3,n)}=\frac{1}{3}\quad\text{and}\quad d_j^{(3,n)}=0\quad\text{for $j\ge3$}.$$
	The difference equation \eqref{eq: DOC-BDF3 difference equation} becomes
	\begin{align*}
		\widehat{\vartheta}_{n-j}^{(3,n)}
		-\frac{7}{11}\widehat{\vartheta}_{n-1-j}^{(3,n-1)}+\frac{2}{11}
		\widehat{\vartheta}_{n-2-j}^{(3,n-2)}=0\quad
		\text{for $n\ge j+2\ge5$.}
	\end{align*}
	Since the associated characteristic equation $\lambda^2-\frac{7}{11}\lambda+\frac{2}{11}=0$
	has two roots $\lambda_{1,2}=(7\pm \imath\sqrt{39})/22$, we have the following explicit formula
	of DOC kernels
	\begin{align*}
		\vartheta_{n-j}^{(3,n)}
		=\frac{11\,\imath}{d_0^{(3,j)}\sqrt{39}}\kbra{\braB{\frac{7-\imath\sqrt{39}}{22}}^{n-j+1}
			-\braB{\frac{7+\imath\sqrt{39}}{22}}^{n-j+1}}
		\quad\text{for $n\ge j+2\ge 5$.}
	\end{align*}
	
	\begin{figure}[htb!]
		\centering
		\includegraphics[height=2in,width=3in]{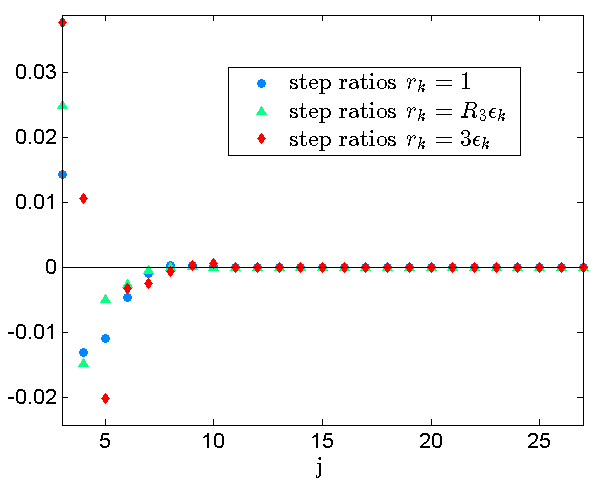}
		\caption{The DOC kernels $\vartheta_{j}^{(3,n)}\,(j\ge2)$
			for different step-ratio patterns.}
		\label{fig: numerical-DOC-kernels}
	\end{figure}
	
	In Figure \ref{fig: numerical-DOC-kernels},
	we fix $n=30$ and list the values of
	DOC kernels $\vartheta_{j}^{(3,n)}$
	for different choices (let $0<\epsilon_{k}<1$ be
	the uniformly distributed random number):
	\begin{itemize}[itemindent=0.4cm]
		\item [(a)] the uniform mesh with $r_k=1$ for $2\leq k\leq 30$;
		\item [(b)] the random mesh with $r_k:=R_3\epsilon_{k}$ for $2\leq k\leq 30$;
		\item [(c)] the random mesh with $r_k:=3 \epsilon_{k}$ for $2\leq k\leq 30$.
	\end{itemize}
	As observed, the DOC kernels $\vartheta_{n-j}^{(3,n)}$ are not always
	positive, but they always decay fast, as predicted by \eqref{eq: DOC-BDF3 decaying estimate}.
	It is interesting to  find that, from the case (c),  the DOC kernels maintain
	the fast damping property as a whole although there are
	many of step ratios greater than $R_3$ (about $15\%$ in the above test).
\end{remark}

\begin{theorem}\label{thm: BDF3 stability}
If $\tau\le 1/(4C_3L_f)$, the BDF3 solution of \eqref{eq: perturbed equation}  satisfies
\begin{align*}
\absb{\tilde{v}^{n}}\le 2\exp(4C_3L_ft_{n-1})\braB{\absb{\tilde{v}^{2}}+5C_3\tau\absb{\partial_{\tau}\tilde{v}^{2}}
	+2C_3\tau\absb{\partial_{\tau}\tilde{v}^{1}}
+2C_3t_{n}\max_{3\le i\le n}\absb{\varepsilon^i}}\quad\text{for $3\le n\le N$.}
\end{align*}
Thus the variable-step BDF3 scheme is stable if $0<r_k<R_3$ for $k\ge2$.
\end{theorem}
\begin{proof}Multiplying both sides of \eqref{eq: DOC perturbed equation} with $2\tau_n\tilde{v}^n$, one
applies Lemma \ref{lem: BDF3 orthogonal formula} to obtain that
\begin{align*}
\absb{\tilde{v}^n}^2-\absb{\tilde{v}^{n-1}}^2
\le&\,-2\tau_n\tilde{v}^n\mathcal{I}_{3}^n[\tilde{v}]
+2\tau_n\tilde{v}^n\sum_{i=3}^n\vartheta_{n-i}^{(3,n)}\kbra{f(t_i,\bar{v}^i)-f(t_i,v^i)+\varepsilon^i}\\
%\le&\,2\tau_n\absb{\tilde{v}^n}\absb{\mathcal{I}_{3}^n[\tilde{v}]}
%+2\tau_n\absb{\tilde{v}^n}\sum_{i=3}^n\absb{\vartheta_{n-i}^{(3,n)}}\absb{f(t_i,\bar{v}^i)-f(t_i,v^i)}
%+2\tau_n\absb{\tilde{v}^n}\max_{3\le i\le n}\absb{\varepsilon^i}\sum_{i=3}^n\absb{\vartheta_{n-i}^{(3,n)}}\\
\le&\,2\tau_n\absb{\tilde{v}^n}\absb{\mathcal{I}_{3}^n[\tilde{v}]}
+2L_f\tau_n\absb{\tilde{v}^n}\sum_{i=3}^n\absb{\vartheta_{n-i}^{(3,n)}}\absb{\tilde{v}^i}
+2C_3\absb{\tilde{v}^n}\tau_n\max_{3\le i\le n}\absb{\varepsilon^i}.
\end{align*}
By summing the above inequality for $n$ from $n=3$ to $m$, we have
\begin{align*}
\absb{\tilde{v}^m}^2\le&\,\absb{\tilde{v}^{2}}^2
+2\sum_{j=3}^m\tau_j\absb{\tilde{v}^j}\absb{\mathcal{I}_{3}^j[\tilde{v}]}
+2L_f\sum_{j=3}^m\tau_j\absb{\tilde{v}^j}\sum_{i=3}^j\absb{\vartheta_{j-i}^{(3,j)}}\absb{\tilde{v}^i}
+2C_3\sum_{j=3}^m\absb{\tilde{v}^j}\tau_j\max_{3\le i\le j}\absb{\varepsilon^i}.
\end{align*}
Taking some integer $m_1$ ($2\le m_1\le m$) such that
$\absb{\tilde{v}^{m_1}}=\max_{2\le k\le m}\absb{\tilde{v}^{k}}$.
Now we take $m:=m_1$ in the above inequality and get
\begin{align}\label{thmproof: BDF3 stability}
\absb{\tilde{v}^{m}}\le \absb{\tilde{v}^{m_1}}\le\absb{\tilde{v}^{2}}
+2\tau\sum_{j=3}^m\absb{\mathcal{I}_{3}^j[\tilde{v}]}
+2C_3L_f\sum_{j=3}^{m}\tau_j\absb{\tilde{v}^j}
+2C_3t_{m}\max_{3\le i\le m}\absb{\varepsilon^i}
\end{align}
for $3\le m\le N$.
It remains to evaluate the error $\sum_{j=3}^m\absb{\mathcal{I}_{3}^j[\tilde{v}]}$ from the starting values.
Taking the index $\rmk=3$ in \eqref{eq: initial effect BDFk-DOC} gives
\begin{align*}
	\mathcal{I}_{3}^n[\tilde{v}]=&\,\partial_{\tau}\tilde{v}^2
	\sum_{i=3}^n\vartheta_{n-i}^{(3,n)}d^{(3,i)}_{i-2}
	+\partial_{\tau}\tilde{v}^1
	\sum_{i=3}^n\vartheta_{n-i}^{(3,n)}d^{(3,i)}_{i-1}\quad\text{for $n\ge3$.}
\end{align*}
Recalling the definition \eqref{eq: BDF3 kernels} of BDF3 kernels with the increasing property \eqref{eq: d0d1d2-increasing}, 
it is easy to check that $d^{(3,3)}_{2}\le d_2(R_3,R_3)\le 2$ and
$\absb{d^{(3,3)}_{1}}+d^{(3,4)}_{2}\le -d_1(R_3,R_3)+d_2(R_3,R_3)\le 5 $.
Thus we apply Lemma \ref{lem: BDF3 orthogonal formula} to get
\begin{align*}
\sum_{j=3}^m\absb{\mathcal{I}_{3}^j[\tilde{v}]}\le&\,
\absb{\partial_{\tau}\tilde{v}^2}
\sum_{j=3}^m\sum_{i=3}^j\absb{\vartheta_{j-i}^{(3,j)}}\absb{d^{(3,i)}_{i-2}}
+\absb{\partial_{\tau}\tilde{v}^1}
\sum_{j=3}^m\sum_{i=3}^j\absb{\vartheta_{j-i}^{(3,j)}}\absb{d^{(3,i)}_{i-1}}\\
=&\,\absb{\partial_{\tau}\tilde{v}^2}
\sum_{i=3}^m\absb{d^{(3,i)}_{i-2}}\sum_{j=i}^m\absb{\vartheta_{j-i}^{(3,j)}}
+\absb{\partial_{\tau}\tilde{v}^1}
\sum_{i=3}^m\absb{d^{(3,i)}_{i-1}}\sum_{j=i}^m\absb{\vartheta_{j-i}^{(3,j)}}\\
\le &\,C_3\absb{\partial_{\tau}\tilde{v}^2}\brab{\absb{d^{(3,3)}_{1}}+d^{(3,4)}_{2}}
+C_3\absb{\partial_{\tau}\tilde{v}^1}d^{(3,3)}_{2}\\
\le &\,5C_3\absb{\partial_{\tau}\tilde{v}^2}+2C_3\absb{\partial_{\tau}\tilde{v}^1}\qquad\text{for $m\ge3$.}
\end{align*}
It follows from \eqref{thmproof: BDF3 stability} that
\begin{align*}
	\absb{\tilde{v}^{n}}\le \absb{\tilde{v}^{2}}
	+5C_3\tau\absb{\partial_{\tau}\tilde{v}^2}+2C_3\tau\absb{\partial_{\tau}\tilde{v}^1}
	+2C_3L_f\sum_{j=3}^{n}\tau_j\absb{\tilde{v}^j}
	+2C_3t_{n}\max_{3\le i\le n}\absb{\varepsilon^i}.
\end{align*}
Assuming that $\tau\le 1/(4C_3L_f)$, one has
\begin{align*}
\absb{\tilde{v}^{n}}\le 2\absb{\tilde{v}^{2}}
+10C_3\tau\absb{\partial_{\tau}\tilde{v}^2}+4C_3\tau\absb{\partial_{\tau}\tilde{v}^1}
+4C_3L_f\sum_{j=3}^{n-1}\tau_j\absb{\tilde{v}^j}
+4C_3t_{n}\max_{3\le i\le n}\absb{\varepsilon^i}
\end{align*}
for $3\le n\le N$. The standard Gr\"{o}nwall inequality completes the proof.
\end{proof}

From the proof of Lemma \ref{lem: BDF3 orthogonal formula},
the uniform boundedness of $\sum_{j=1}^{n}\absb{\vartheta_{n-j}^{(3,n)}}$
and $\sum_{j=i}^{n}\absb{\vartheta_{j-i}^{(3,j)}}$
does not require all DOC kernels to decay rapidly.
It always allows a finite number of bounded DOC kernels.
As seen from the numerical tests in the next section, 
the imposed step-ratio condition $0<r_k<R_3$
in Theorem \ref{thm: BDF3 stability} is sufficient, but far from necessary.
A practical step-ratio constraint for stability is that
\begin{center}
\emph{almost all step ratios satisfy $0<r_k<R_3$ for $2\le k\le N$}.
\end{center}
%This is a proper understanding of the constraint $0<r_k<R_3$ in
%Lemma \ref{lem: BDF3 orthogonal formula}
%and Theorem \ref{thm: BDF3 stability}.

\section{Numerical tests}
\setcounter{equation}{0}

Consider an ODE model $v'=2v-3e^{-t}$ for $0<t\le 1$ with a smooth solution $v=\exp(-t)$ in our numerical tests.
We run the variable-step BDF2 and BDF3 schemes in two scenarios:
\begin{enumerate}
	\item[(a)] The graded meshes $t_k=(k/N)^\gamma$ for $0\le k\le N$ with grading parameters $\gamma>1$.
	The maximum step-ratio is $r_{\max}=r_2=2^{\gamma}-1$ and 
	$\tau/\tau_1=\tau_N/\tau_1=N^{\gamma}-(N-1)^{\gamma}\approx \gamma N^{\gamma-1}$.
	\item[(b)] The random meshes with time-steps $\tau_k=\epsilon_k/\sum_{k=1}^N\epsilon_k$, where $\epsilon_k\in (0,1)$ are uniformly distributed random numbers.
\end{enumerate}
To start the multi-step schemes, a two-stage third-order singly diagonally 
implicit Runge-Kutta method is employed. 
We record the maximum error $e(N):=\max_{1\leq{n}\leq{N}}|v(t_n)-v^n|$
in each run and compute the convergence order by
$$\text{Order}\approx\frac{\log\bra{e(N)/e(2N)}}{\log\bra{\tau(N)/\tau(2N)}}$$
where $\tau(N)$ denotes the maximum time-step size for total $N$ subintervals.

%%%%%%%%%%%%%%%%%%%%%%%%%%%%%%%%%%%%%%%%%%%%%%%%%%%%%%%%%%%%%%%%%%%%%%%%%%%%

%%%%%%%%%%%%%%%%%%%%%%%%%%%%%%%%%%%%%%%%%%%%%%%%%%%%%%%%%%%%%%%%%%%%%%%%%%%
\begin{table}[htb!]
	\begin{center}
		\caption{Numerical results of BDF2 method on graded time meshes}\label{table:BDF2 test on graded meshes} \vspace*{0.3pt}
		\def\temptablewidth{1.0\textwidth}
		{\rule{\temptablewidth}{0.5pt}}
		\begin{tabular*}{\temptablewidth}{@{\extracolsep{\fill}}cccccccccc}
			\multirow{2}{*}{$N$} &\multirow{2}{*}{$\dfrac{\tau}{\tau_{1}}$} &\multicolumn{2}{c}{$\gamma=2,\,\,r_{\max}=3$} &\multirow{2}{*}{$\dfrac{\tau}{\tau_{1}}$} &\multicolumn{2}{c}{$\gamma=3,\,\,r_{\max}=7$} &\multirow{2}{*}{$\dfrac{\tau}{\tau_{1}}$}&\multicolumn{2}{c}{$\gamma=4,\,\,r_{\max}=15$} \\
			\cline{3-4}          \cline{6-7}         \cline{9-10}
			&          &$e(N)$   &Order &         &$e(N)$   &Order &          &$e(N)$    &Order\\
			\midrule
			40     &7.9e+01  &5.28e-04   &$-$   &4.7e+03   &8.77e-04     &$-$   &2.5e+05  &1.35e-03  &$-$\\
			80     &1.6e+02  &1.34e-04	 &1.97	&1.9e+04   &2.25e-04	 &1.96  &2.0e+06  &3.49e-04	 &1.95\\
			160    &3.2e+02  &3.39e-05	 &1.99  &7.6e+04   &5.72e-05	 &1.98  &1.6e+07  &8.91e-05	 &1.97\\
			320    &6.4e+02  &8.52e-06	 &1.99  &3.1e+05   &1.44e-05	 &1.99  &1.3e+08  &2.25e-05	 &1.98\\
			640    &1.3e+03  &2.14e-06	 &2.00	&1.2e+06   &3.61e-06	 &1.99  &1.1e+09  &5.66e-06	 &1.99\\
			1280   &2.6e+03  &5.34e-07	 &2.00	&4.9e+06   &9.06e-07	 &2.00  &8.4e+09  &1.42e-06	 &2.00\\
			%			\midrule
			%			\multicolumn{3}{l}{$\min\{1\alpha, \gamma\sigma\}$}   &1.20 & & &1.80 & & &1.80\\
		\end{tabular*}
		{\rule{\temptablewidth}{0.5pt}}
	\end{center}
\end{table}	
%%%%%%%%%%%%%%%%%%%%%%%%%%%%%%%%%%%%%%%%%%%%%%%%%%%%%%%%%%%%%%%%%%%%%%%%%%%%

%%%%%%%%%%%%%%%%%%%%%%%%%%%%%%%%%%%%%%%%%%%%%%%%%%%%%%%%%%%%%%%%%%%%%%%%%%%
\begin{table}[htb!]
	\begin{center}
		\caption{Numerical results of BDF3 method on graded time meshes}\label{table:BDF3 test on graded meshes} \vspace*{0.3pt}
		\def\temptablewidth{1.0\textwidth}
		{\rule{\temptablewidth}{0.5pt}}
		\begin{tabular*}{\temptablewidth}{@{\extracolsep{\fill}}cccccccccc}
			\multirow{2}{*}{$N$} &\multirow{2}{*}{$\dfrac{\tau}{\tau_{1}}$} &\multicolumn{2}{c}{$\gamma=2,\,\,r_{\max}=3$} &\multirow{2}{*}{$\dfrac{\tau}{\tau_{1}}$} &\multicolumn{2}{c}{$\gamma=3,\,\,r_{\max}=7$} &\multirow{2}{*}{$\dfrac{\tau}{\tau_{1}}$}&\multicolumn{2}{c}{$\gamma=4,\,\,r_{\max}=15$} \\
			\cline{3-4}          \cline{6-7}         \cline{9-10}
			&          &$e(N)$   &Order &         &$e(N)$   &Order &          &$e(N)$    &Order\\
			\midrule
			40     &7.9e+01  &1.27e-05  &$-$    &4.7e+03   &2.94e-05 &$-$   &2.5e+05  &5.73e-05  &$-$\\
			80     &1.6e+02  &1.65e-06	 &2.95	&1.9e+04   &3.91e-06 &2.91  &2.0e+06  &7.85e-06  &2.87\\
			160    &3.2e+02  &2.10e-07	 &2.97  &7.6e+04   &5.05e-07 &2.96  &1.6e+07  &1.03e-06  &2.94\\
			320    &6.4e+02  &2.65e-08	 &2.99  &3.1e+05   &6.41e-08 &2.98  &1.3e+08  &1.31e-07  &2.97\\
			640    &1.3e+03  &3.32e-09	 &2.99	&1.2e+06   &8.07e-09 &2.99  &1.1e+09  &1.66e-08  &2.98\\
			1280   &2.6e+03  &4.16e-10	 &3.00	&4.9e+06   &1.01e-09 &2.99  &8.4e+09  &2.08e-09  &2.99\\
			%			\midrule
			%			\multicolumn{3}{l}{$\min\{1\alpha, \gamma\sigma\}$}   &1.20 & & &1.80 & & &1.80\\
		\end{tabular*}
		{\rule{\temptablewidth}{0.5pt}}
	\end{center}
\end{table}	
%%%%%%%%%%%%%%%%%%%%%%%%%%%%%%%%%%%%%%%%%%%%%%%%%%%%%%%%%%%%%%%%%%%%%%%%%%%%

Tables \ref{table:BDF2 test on graded meshes}-\ref{table:BDF3 test on graded meshes} 
list the numerical results of 
the BDF2 and BDF3 methods on graded meshes with three grading parameters, respectively.
We see that,
although the ratio $\tau/\tau_1$ increases as fast as $O(N^{\gamma-1})$, 
the numerical solutions remains stable and convergent with full accuracy.

%%%%%%%%%%%%%%%%%%%%%%%%%%%%%%%%%%%%%%%%%%%%%%%%%%%%%%%%%%%%%%%%%%
\begin{table}[htb!]
	\begin{center}
		\caption{Numerical results of BDF2 method on random time meshes}
		\label{table:BDF2 test on random meshes}\vspace*{0.3pt}
		\def\temptablewidth{0.75\textwidth}
		{\rule{\temptablewidth}{0.5pt}}
		\begin{tabular*}{\temptablewidth}{@{\extracolsep{\fill}}cccccc}
			$N$   &$\tau(N)$      &$e(N)$     &Order  &$r_{\max}$  &$\tau/\tau_{1}$\\
			\midrule
			40    &4.49e-02	 &5.86e-04	 &--	 &28.30   &1.22\\
			80    &2.40e-02	 &1.43e-04	 &2.03	 &91.41   &1.06\\
			160   &1.18e-02	 &4.09e-05	 &1.81	 &32.54   &19.88\\
			320	  &6.18e-03	 &1.03e-05	 &1.99	 &418.41  &1.21\\	
			640	&3.01e-03	 &2.51e-06	 &2.04	 &604.02  &2.66\\
			1280 &1.55e-03	 &6.63e-07	 &1.92	 &1963.80 &1.02
		\end{tabular*}
		{\rule{\temptablewidth}{0.5pt}}
	\end{center}
\end{table}		
%%%%%%%%%%%%%%%%%%%%%%%%%%%%%%%%%%%%%%%%%%%%%%%%%%%%%%%%%%%%%%%%%%

%%%%%%%%%%%%%%%%%%%%%%%%%%%%%%%%%%%%%%%%%%%%%%%%%%%%%%%%%%%%%%%%%%
\begin{table}[htb!]
	\begin{center}
		\caption{Numerical results of BDF3 method on random time meshes}
		\label{table:BDF3 test on random meshes}\vspace*{0.3pt}
		\def\temptablewidth{0.75\textwidth}
		{\rule{\temptablewidth}{0.5pt}}
		\begin{tabular*}{\temptablewidth}{@{\extracolsep{\fill}}ccccccc}
			$N$   &$\tau(N)$      &$e(N)$     &Order  &$r_{\max}$   &$N_1$ &$\tau/\tau_{1}$\\
			\midrule
			40    &4.60e-02	 &8.20e-06	 &--	 &7.58     &9     &1.57\\
			80    &2.46e-02	 &1.18e-06	 &2.79	 &361.49   &18    &1.93\\
			160   &1.22e-02	 &1.83e-07	 &2.69	 &1682.18  &35    &4.51\\
			320	  &6.11e-03	 &2.40e-08	 &2.93	 &79.90    &60    &2.91\\	
			640	  &3.27e-03	 &3.40e-09	 &2.82	 &5765.00  &146   &1.20\\
			1280  &1.56e-03	 &4.16e-10	 &3.03	 &9677.92  &250   &6.64    
		\end{tabular*}
		{\rule{\temptablewidth}{0.5pt}}
	\end{center}
\end{table}		
%%%%%%%%%%%%%%%%%%%%%%%%%%%%%%%%%%%%%%%%%%%%%%%%%%%%%%%%%%%%%%%%%%

Tables \ref{table:BDF2 test on random meshes}-\ref{table:BDF3 test on random meshes} 
record the numerical results on random time meshes.
Table \ref{table:BDF2 test on random meshes} suggests that
the variable-step BDF2 method is robust with respect to the step-ratios $r_k$,
and supports the theoretical findings in Theorem 2.1.
Reminding the step-ratio restriction $0<r_k<2.533$ in Theorem 3.1 for the BDF3 method,
we also record the number (denoted by $N_1$ in
Table \ref{table:BDF3 test on random meshes}) of time levels with $r_k\ge 2.533$. 
It is seen that the variable-step BDF3 method is stable and third-order convergent 
on random meshes, even if there are many of (about $20\%$ in our tests) 
step-ratios do not meet our theoretical condition.
Nonetheless, it remains mysterious to us.

\section{Concluding remarks}
\setcounter{equation}{0}

The stability of BDF2 and BDF3 methods with unequal time-steps is
investigated by a new theoretical framework using the discrete orthogonal convolution kernels.
Thanks to the global analysis, that is, 
the present technique formulates the current BDF solution as a convolution summation of 
all previous information with DOC kernels as 
the convolutional weights, see the form \eqref{eq: DOC perturbed equation},
we improved the setp-ratio constraints for the stability. 
It is to mention that, the global decaying estimates of DOC kernels
are also critical in the numerical analysis of BDF methods for parabolic equations,
cf. \cite{LiaoJiZhang:2020pfc,LiaoJiWangZhang:2021,LiaoZhang:2021}. 
The stability and convergence theory of variable-step BDF3 scheme for these stiff problems
will be addressed in forthcoming reports.

For the high-order BDF-$\rmk$ ($\rmk=4,5,6$) time-stepping methods,
we can also define the associated DOC kernels $\vartheta_{n-j}^{(\rmk,n)}$ via
 the recursive procedure \eqref{eq: BDFk-DOC procedure}.
From the proof of Theorem \ref{thm: BDF3 stability},
one needs a result similar to Lemma \ref{lem: BDF3 orthogonal formula}
under certain step-ratio condition.
Actually, the uniform boundedness (there exists a constant $C_{\rmk}$
independent of the time-level indexes $n$)
of the absolute summations of DOC kernels
\begin{align*}
	\sum_{j=\rmk}^{n}\absb{\vartheta_{n-j}^{(\rmk,n)}}\le C_{\rmk}\quad\text{and}\quad
	\sum_{j=i}^{n}\absb{\vartheta_{j-i}^{(\rmk,j)}}\le C_{\rmk}\quad\text{for any $n\ge i\ge\rmk$,}
\end{align*}
is sufficient to the stability of BDF-$\rmk$ schemes for initial value problems.
However, this issue seems to be theoretically challenging for $\rmk\ge4$
and remains open to us.

\section*{Acknowledgements}
The authors would like to thank Dr. Ji Bingquan for his valuable
discussions.

\appendix
\section{Two auxiliary functions for the BDF3 method}
\label{append: three auxiliary functions}
\setcounter{equation}{0}

By using the definitions \eqref{funs: d0}-\eqref{funs: d2}, it is easy to check that
\begin{align}
	&\frac{\partial }{\partial x}\abs{d_{\nu}(x,y)}
>0\quad\text{and}\quad \frac{\partial}{\partial y}\abs{d_{\nu}(x,y)}
>0\quad\text{for $\nu=0,1,2$ and $x,y>0.$}\label{eq: d0d1d2-increasing}
\end{align}
That is, the functions $d_0(x,y)$, $-d_1(x,y)$ and $d_2(x,y)$ are increasing with respect to $x,y>0$.
We consider two auxiliary functions for the analysis of variable-step BDF3 method,
\begin{align}
\alpha(x,y):=&\,-\frac{d_1(x,y)}{d_0(x,y)}=\frac{x\left(x^2y^2+4 xy^2+3y^2+2 xy +3y+1\right)}
{(y+1) \left(3 x^2 y+4xy+2x+y+1\right)},
\label{eq: alpha auxiliary function}\\
\beta(x,y):=&\,\frac{d_2(x,y)}{d_0(x,y)}=\frac{x (x+1)^2y^2 }{(y+1)\left(3 x^2 y+4xy+2x+y+1\right)}.
\label{eq: beta auxiliary function}
\end{align}
Their properties will be examined  with respect to
two independent variables $x,y>0$ due to
the facts $\alpha_m=\alpha(r_{m},r_{m-1})$ and $\beta_m=\beta(r_{m},r_{m-1})$
for $m\ge2$ according to \eqref{eq: DOC-BDF3 coefficients}.

The above functions $\alpha$ and $\beta$ are strictly increasing with respect to $x,y>0$.
Actually, by simple but tedious calculations for \eqref{eq: alpha auxiliary function}
-\eqref{eq: beta auxiliary function}, it is not difficult to check that
\begin{align}
&\frac{\partial \alpha}{\partial x}
=\frac{(x+1)^2 \left(3 x^2+2 x+3\right) y^3+2 \left(2 x^3+5 x^2+6 x+3\right) y^2+(x+2)^2 y+1}
{(y+1) \left(3 x^2 y+4 xy+2x+y+1\right)^2},\label{eq: alpha x-increasing}\\
&\frac{\partial \alpha}{\partial y}=\frac{x (x+1)^2 (x y+y+1) (3xy+3y+1)}
{(y+1)^2 \left(3 x^2 y+ 4x y+2x+y+1\right)^2},
\label{eq: alpha y-increasing}
\end{align}
and
\begin{align}
&\frac{\partial \beta}{\partial x}
=\frac{(x+1) y^2 \left(3 x^3 y+ 5x^2 y+4x^2+3 xy+ 3x+y+1\right)}
{(y+1) \left(3 x^2 y+4 xy+2x+y+1\right)^2},\label{eq: beta x-increasing}\\
&\frac{\partial \beta}{\partial y}
=\frac{x (x+1)^2 y \left(3 x^2 y+6 x y+4 x+2 y+2\right)}
{(y+1)^2 \left(3 x^2 y+4x y+2x+y+1\right)^2}\quad\text{for $x,y>0$.}
\label{eq: beta y-increasing}
\end{align}

We have the following results. Some of proofs are technically complex and 
the mathematical derivations have been checked carefully by a symbolic calculation software.
\begin{lemma}\label{lemma: alpha beta 1}
It holds that
$$0<\beta(x,y)<\alpha(x,y)<1+\beta(x,y)\quad\text{for $x,y>0$}.$$
\end{lemma}
\begin{proof}Simple calculations lead to
\begin{align*}
\alpha(x,y)-\beta(x,y)=&\,\frac{x (2xy+2y+1)}{3 x^2 y+4x y+2x+y+1},\\
1-\alpha(x,y)+\beta(x,y)=&\,\frac{(x+1) (x y+y+1)}{3 x^2 y+4x y+2x+y+1}.
\end{align*}
It completes the proof.
\end{proof}
\begin{lemma}\label{lemma: beta less than 1}
It holds that
$$0<\beta(x,y)<1\quad\text{for $0<x,y<\hat{R}_3\approx3.4405$,}$$
where $\hat{R}_3$ is the unique positive root of
$\hat{R}_3^4-2 \hat{R}_3^3-4 \hat{R}_3^2-3 \hat{R}_3-1=0$.
\end{lemma}
\begin{proof}According to
\eqref{eq: beta x-increasing}-\eqref{eq: beta y-increasing},
it holds that
\begin{align*}
\beta(x,y)-1<\beta(\hat{R}_3,\hat{R}_3)-1
=\frac{\hat{R}_3^4-2 \hat{R}_3^3-4 \hat{R}_3^2-3 \hat{R}_3-1}
{3 \hat{R}_3^3+4 \hat{R}_3^2+3 \hat{R}_3+1}=0\quad\text{for $0<x,y<\hat{R}_3$.}
\end{align*}
It completes the proof.
\end{proof}

\begin{lemma}\label{lemma: elliptic function}
For $0< x,y< R_3\approx2.553$, it holds that
$$2\alpha^2(x,y)+3\beta^2(x,y)-4\alpha(x,y)\beta(x,y)-2\alpha(x,y)+2\beta(x,y)<0.$$
\end{lemma}

\begin{proof}
Consider an auxiliary function, see Figure \ref{fig: surface g},
\begin{align*}
g(x,y):=&\,\frac{1}{x (x+1)}\kbraB{2\alpha(x,y)^2+3\beta(x,y)^2-4\alpha(x,y)\beta(x,y)
-2\alpha(x,y)+2\beta(x,y)}\\
=&\,\frac{\left(3 x^2 y+4x y+2x+y+1\right)}{\gamma(x,y)}
\Big[(x+1)^2 \left(x^2+x-4\right) y^4-2 \left(4 x^2+11 x+7\right) y^3\Big.\\
&\,\hspace{3.5cm}\Big.-2 \left(2 x^2+10 x+9\right) y^2-2 (3 x+5) y-2\Big]\qquad\text{for $0\le x,y< R_3$,}
\end{align*}
where
$$\gamma(x,y):=(y+1)^2 \left(3 x^2 y+4x y+2x+y+1\right)^3.$$

\begin{figure}[htb!]
\centering
\includegraphics[height=2in,width=3in]{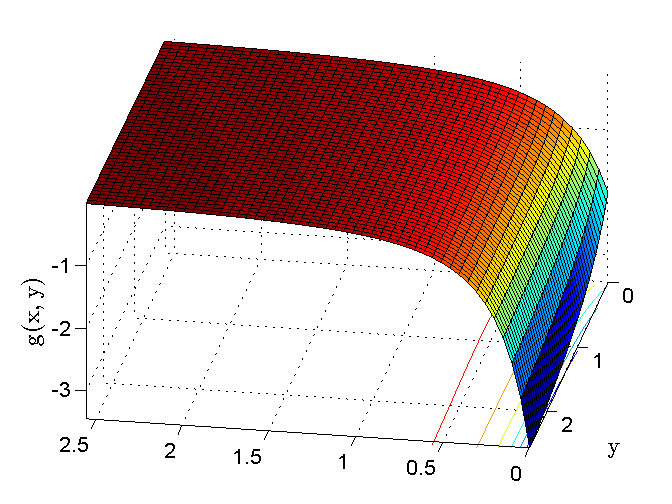}
\caption{The surface of $g$ on $[0,R_3]^2$.}
\label{fig: surface g}
\end{figure}

Obviously, it is sufficient to show that $g(x,y)<0$ for $0< x,y< R_3$.
Simple but tedious calculations yield
\begin{align*}
\frac{\partial g}{\partial x}=&\,\frac{1}{\gamma(x,y)}
\Big[2 \left(27 x^2+88 x+63\right) y^2+2 \left(12 x^3+90 x^2+158 x+79\right) y^3\Big.\\
&\,\hspace{1cm}\Big.+\left(4 x^4+58 x^3+207 x^2+252 x+99\right) y^4+(25-x) (x+1)^3 y^5+(36 x+50) y+8\Big].
\end{align*}
Obviously, we have $\frac{\partial g}{\partial x}>0$ such that
$g$ is increasing with respect to $x$ for $0<x,y<R_3$.
Thus it has no extreme points over the open square $(0,R_3)^2$.
It remains to consider the maximum value along the four sides of $[0,R_3]^2$:
\begin{itemize}
  \item [(i)]Along the side $y=0$, we have $g(x,0)=-\frac{2}{(2 x+1)^2}<0$ for $0< x< R_3$..
  \item [(ii)]Along the side $y=R_3$, $g(x,R_3)<g(R_3,R_3)\approx-0.0000277<0$ for $0< x< R_3$.
  \item [(iii)] Along the side $x=0$, we have $g(0,y)=-\frac{4 y+2}{y+1}<0$ for $0< y< R_3$..
  \item [(iv)] Along the side $x=R_3$, one has
  \begin{align*}
g(R_3,y)\approx&\,\frac{0.06763y^4-0.12922 y^3
-0.100507 y^2-0.0267488 y-0.002113}{(y+1)^2 (y+0.19847)^2}
\quad\text{for $0< y< R_3$.}
\end{align*}
It has a unique minimum point at $y\approx0.103652$, while
$$g(R_3,y)<\max\{g(R_3,0),g(R_3,R_3)\}=g(R_3,R_3)<0\quad\text{for $0< y< R_3$.}$$
\end{itemize}
In summary, we have $g(x,y)<0$ for $0< x,y< R_3$. It completes the proof.
\end{proof}

%\section{Tests for the BDF4 method}
%\begin{align*}
%b_0(x,y,z):=&\,1+\frac{x}{1+x}+\frac{xy}{1+y+xy}+\frac{xyz}{1+z+yz+xyz}\\
%b_1(x,y,z):=&\,-1-x-\frac{xy(1+x)}{1+y}
%-xyz\frac{1+x}{1+y}\frac{1+y+xy}{1+z+yz}\\
%b_2(x,y,z):=&\,\frac{x^2}{1+x}+xy^2+xy^2(1+y+xy)\frac{z}{1+z}\\
%b_3(x,y,z):=&\,-x^2y^3\frac{1+z+yz+xyz}{1+y+xy}
%\frac{1+x}{1+y}\\
%b_4(x,y,z):=&\,\frac{x^2y^3z^4}{1+z+yz+xyz}\frac{1+x}{1+z}\frac{1+y+xy}{1+z+yz}
%\qquad\text{for $x, y\ge0$.}
%\end{align*}

\end{document}